\documentclass{amsart}
\usepackage{palatino,hhline, bbm}
\usepackage{arydshln} 
\usepackage{amsthm, amsmath}
\usepackage{amssymb} 
\usepackage{amscd}
\usepackage{mathrsfs}
\renewcommand{\mathcal}{\mathscr}
\usepackage[hmargin=3cm, vmargin=3cm]{geometry}
\usepackage[colorlinks,bookmarksopen,bookmarksnumbered, linkcolor=blue,citecolor=red,urlcolor=blue]{hyperref}
\usepackage[all]{xy}

\DeclareMathAlphabet{\mathpzc}{OT1}{pzc}{m}{it}

\theoremstyle{plain}

\newtheorem*{thm*}{Theorem}

\newtheorem{lem}{Lemma}
\newtheorem{prop}{Proposition}

\theoremstyle{remark}

\newcommand\pr{\noindent\textit{Proof} : }

\newcommand\rond{\kern 1pt{\scriptstyle\circ}\kern 1pt}

\def\lr#1{\langle {#1} \rangle}

\newcommand\Ext{\operatorname{Ext}}

\newcommand\Hom{\operatorname{Hom}}
\newcommand\im{\operatorname{Im}}

\newcommand\Pic{\operatorname{Pic}}

\newcommand\Tr{\operatorname{Tr}}
\newcommand\pp{\scriptscriptstyle\bullet}

\newcommand\Z{\mathbb{Z}}

\newcommand\C{\mathbb{C}}
\renewcommand\P{\mathbb{P}}

\newcommand\G{\mathbb{G}}
\renewcommand\O{\mathcal{O}}

\newcommand\iso{\vbox{\hbox to .8cm{\hfill{$\scriptstyle\sim$}\hfill}
\nointerlineskip\hbox to .8cm{{\hfill$\longrightarrow $\hfill}} }}
\newcommand\bir{\vbox{\hbox to .8cm{\hfill{$\scriptstyle\sim$}\hfill}
\nointerlineskip\hbox to .8cm{{\hfill$\dasharrow $\hfill}} }}
\newcommand\abs[1]{\lvert {#1}\rvert}

\begin{document}
\title{Vector bundles on Fano threefolds and K3 surfaces}
\author[Arnaud Beauville]{Arnaud Beauville}
\address{Universit\'e C\^ote d'Azur\\
CNRS -- Laboratoire J.-A. Dieudonn\'e\\
Parc Valrose\\
F-06108 Nice cedex 2, France}
\email{arnaud.beauville@unice.fr}
 
\begin{abstract}
Let $X$ be a Fano threefold, and let $S\subset X$ be a K3 surface. Any moduli space $\mathcal{M}_{S}$  of simple vector bundles on $S$ carries a holomorphic symplectic structure. Following an idea of Tyurin, we show that in some cases, those vector bundles which come from $X$ form a Lagrangian subvariety of $\mathcal{M}_{S}$. We illustrate this with a number of concrete examples.

\end{abstract}
\maketitle 

\section*{Introduction}
Let $X$ be a Fano threefold, and let $S$ be a smooth surface in the anticanonical system of $X$, so that $S$ is a K3 surface. Suppose we have a nice moduli space of vector bundles $\mathcal{M}_X$ on $X$, such that their restriction to $S$ belongs to a moduli space $\mathcal{M}_S$. Under mild hypotheses $\mathcal{M}_S$ has a natural (holomorphic) symplectic structure. What can we say of the restriction map $\operatorname{res} :\mathcal{M}_X\rightarrow \mathcal{M}_S$, in particular with respect to this symplectic structure?

In a 1990 preprint \cite{T}, Tyurin made a remarkable observation: if $H^2(X,\mathcal{E}nd(E))=0$ for all $E$ in $\mathcal{M}_{X}$, $\operatorname{res} $ is a local isomorphism to a Lagrangian subvariety of $\mathcal{M}_S$. The proof is quite simple (see \S 1). However Tyurin does not give any example  where this result can be applied. 
We will show that under appropriate hypotheses, the Serre construction provides  such a situation (in rank 2). This will give us a number of examples  of Lagrangian subvarieties, in particular inside the O'Grady hyperk\"ahler manifold $\mathrm{OG}_{10}$.

\medskip	
\hskip0.5cm \Small{I am grateful to Daniele Faenzi for pointing out the reference \cite{B-F}, which has allowed me to extend and simplify the results of \S 7.}

\section{Tyurin's theorem}
Throughout the paper, we will denote by $X$ a Fano threefold (over $\C$), and by $S$ a smooth surface in the anticanonical system $\abs{K_X^{-1}}$; thus $S$ is a K3 surface. Recall that the moduli space  of simple vector bundles on $X$ or $S$ exists as an algebraic space \cite{A-K} (equivalently, as an analytic Moishezon space). 
\begin{thm*}[Tyurin] Let $\mathcal{M}_X$ be a component of the moduli space  of simple vector bundles on $X$.
Assume $H^2(X,\mathcal{E}nd(E))=0$ for every $E$ in $\mathcal{M}_{X}$. Then:

$1)$ $\mathcal{M}_X$ is smooth; for each $E$ in $\mathcal{M}_{X}$, the vector bundle $E_{|S}$ is simple. 

\smallskip	
\hspace{1cm}\parbox{14cm}{\indent\emph{Let $\mathcal{M}_{S}$ be the component of the moduli space of simple vector bundles on $S$ containing the   vector bundles $E_{|S}$. By \cite{M1} $\mathcal{M}_S$ is smooth and carries a symplectic structure.}}

\smallskip	
 2)  The restriction map $\operatorname{res} :\mathcal{M}_{X}\rightarrow \mathcal{M}_{S}$ induces a local isomorphism of $\mathcal{M}_{X}$ into a Lagrangian subvariety of $\mathcal{M}_{S}$.

\end{thm*}
\pr  
Let $E$ be a vector bundle in $\mathcal{M}_{X}$, and let $E_{S}$ be its restriction to $S$. The condition \break $H^2(X,\mathcal{E}nd(E))=0$ implies that $\mathcal{M}_X$ is smooth at $[E]$. Let $s$ be a section of $K_{X}^{-1}$ defining $S$. 
Tensoring the exact sequence $0\rightarrow K_{X}\xrightarrow{\ s\ }\O_{X}\rightarrow \O_S\rightarrow 0$ with $\mathcal{E}nd(E)$ gives an exact sequence
\[0\rightarrow \mathcal{E}nd(E)\otimes K_{X}\rightarrow \mathcal{E}nd(E)\rightarrow \mathcal{E}nd(E_S)\rightarrow 0\,.\]
Consider the associated long exact sequence. Since $H^2(X,\mathcal{E}nd(E))$ and its dual $H^1(X,\mathcal{E}nd(E)\otimes K_X)$ are zero,  we get that the restriction map $H^0(S,\mathcal{E}nd(E))\rightarrow H^0(S,\mathcal{E}nd(E_S))$ is an isomorphism, hence $E_S$ is simple. We also get
 an exact sequence
\[0\rightarrow H^1(X,\mathcal{E}nd(E))\rightarrow H^1(S,\mathcal{E}nd(E_S))\rightarrow H^2(X,\mathcal{E}nd(E)\otimes K_{X})\rightarrow 0\, .\]The first space is isomorphic to the tangent space $T_{E}(\mathcal{M}_X)$, and the second to $T_{E_S}(\mathcal{M}_S)$; 
this shows that  $\operatorname{res} $ is an immersion. Moreover, we get
$\dim T_{E}(\mathcal{M}_X) =  \frac{1}{2}\dim T_{E_S}(\mathcal{M}_S)$ by Serre duality. From the commutative diagram 
\[\xymatrix@M=5pt{H^1(X,\mathcal{E}nd(E)))\otimes H^1(X,\mathcal{E}nd(E)))\ar[r] \ar[d]&H^1(S,\mathcal{E}nd(E_S)))\otimes H^1(S,\mathcal{E}nd(E_S)))\ar[d]\\
H^2(X,\mathcal{E}nd(E))) \ar[r]\ar[d]^{\Tr}& H^2(S,\mathcal{E}nd(E_S)))\ar[d]^{\Tr}\\
0=H^2(X,\O_X)\ar[r]& H^2(S,\O_S)\cong \C}\]
we deduce that the image of $T_{E}(\mathcal{M}_X)$ in $  T_{E_S}(\mathcal{M}_S)$ is isotropic, hence Lagrangian.\qed

\medskip	
\noindent\emph{Remarks}$.-$
1) Without the vanishing hypothesis, the proof still shows  that the image of $T_E(\mathcal{M}_X)$ in $T_{E_S}(\mathcal{M}_S)$ is  Lagrangian; however, if $E_S$ is simple, the tangent map $T(\operatorname{res} ):T_E(\mathcal{M}_X)\rightarrow T_{E_S}(\mathcal{M}_S)$ is injective if and only if $H^2(\mathcal{E}nd(E))=0$.

2) Let $r$ be the rank of the bundles $E$ in $\mathcal{M}_X$, and let $\Delta :=2rc_2-(r-1)c_1^2$ be their discriminant. Under the hypotheses of the theorem, we have by Riemann-Roch 
\begin{equation}\label{RR}
\dim \mathcal{M}_X=h^1(\mathcal{E}nd(E))=1-\chi (E)=\dfrac{1}{2}(-K_X\cdot \Delta _E)+1-r^2 \, .
\end{equation}

\medskip	
\section{The Serre construction}
Let $X$ be a Fano threefold, $C\subset X$ a smooth curve (or, more generally, a locally complete intersection curve). We assume that there exists an ample line bundle $L$ on $X$ such that $(K_{X}\otimes L)_{|C}\cong K_{C}$. There is a unique rank 2 vector bundle $E$ on $X$ and an extension 
\begin{equation}\label{ext}
0\rightarrow \O_X\rightarrow E\rightarrow \mathcal{I}_CL\rightarrow 0\,.
\end{equation}(see e.g. \cite[Remark 1.1.1]{Har}). 
Restricting to $C$ gives an isomorphism  of $E_{|C}$ onto the normal bundle $ N_{C}$ of $C$ in $X$.
\begin{prop}\label{main}
Assume that $H^1(C,N_{C})=0$, and that the restriction map $H^0(X,K_X\otimes L)\rightarrow H^0(C,K_C)$ is surjective. 
Then  $H^2(X,\mathcal{E}nd(E))=0$.
\end{prop}
\pr Tensoring (\ref{ext}) with $E^*$ gives an exact sequence
\begin{equation}\label{end}
0\rightarrow E^*\rightarrow  \mathcal{E}nd(E) \rightarrow \mathcal{I}_{C}E \rightarrow  0\, .
\end{equation}
To prove $H^2(\mathcal{E}nd(E))=0$, we will prove that both $H^2(E^*)$ and $H^2(\mathcal{I}_CE)$ are zero.

1) We have $H^2(E^*)\cong H^1(E\otimes K_{X})^*$ by Serre duality. The exact sequence (\ref{ext}) gives an isomorphism $H^1(E\otimes K_{X})\iso \allowbreak H^1(\mathcal{I}_{C} K_{X}\otimes L)$. Since $H^1(K_{X}\otimes L)=0$ by Kodaira vanishing, we get from the exact sequence $0\rightarrow \mathcal{I}_{C}\rightarrow \O_X\rightarrow \O_{C}\rightarrow 0$ an isomorphism of $H^1(\mathcal{I}_{C} K_{X}\otimes L)$ onto the cokernel of the restriction map  $H^0(X,K_X\otimes L)\rightarrow H^0(C,K_C)$, which is zero by our hypothesis.

2)  $H^2(\mathcal{I}_{C}E)$ fits into an exact sequence $H^1(E_{|C})\rightarrow H^2(\mathcal{I}_{C}E)\rightarrow H^2(E)$.  We have $H^1(E_{|C})=\allowbreak H^1(N_{C})=0$ by hypothesis. Using again (\ref{ext}) we find $H^2(E)\cong H^2(\mathcal{I}_{C}L)$. Since $H^i(L)=0$ for $i=1,2$, this is isomorphic to $H^1(L_{|C})=H^1(K_{C}\otimes K_{X}^{-1}{}_{|C})$, which is zero since $K_{X}^{-1}$ is ample.\qed

\begin{prop}\label{deform}
Assume that the restriction map $H^0(X,L)\rightarrow H^0(C,L_{|C})$ is surjective. Then the vector bundles $E$ obtained by the above construction form an open subset of the moduli space $\mathcal{M}_X$.
\end{prop}
\pr Let $E$ be the vector bundle in (\ref{ext}), and $s$ a section of $E$ vanishing along $C$. The tangent space to the deformation space $\operatorname{Def}(E,s) $ of the pair $(E,s)$ is  $H^1(X,K^{\pp})$, where $K^{\pp}$ is the complex $\mathcal{E}nd(E) \xrightarrow{\ e_s\ } E$, with $e_s(u)=u(s)$. 
The tangent space to $\operatorname{Def}_{E} $ at $E$ is $H^1(X,\mathcal{E}nd(E))$, and the tangent map to the forgetful map $\operatorname{Def}(E,s)\rightarrow  \operatorname{Def}(E)$ is $H^1(p)$, where $p:K^{\pp}\rightarrow \mathcal{E}nd(E))$ is the projection; it suffices to prove that $H^1(p)$ is surjective. Now $p$ is surjective with kernel $E[-1]$, so it suffices to prove that $H^1(E)$ is zero. From the exact sequence (\ref{ext}) this is equivalent to $H^1(\mathcal{I}_CL)=0$, which in turn is equivalent to our hypothesis. \qed

\medskip	
Let $\mathcal{C}_X$ be the component of the Hilbert scheme parametrizing the curves $C\subset X$. The Serre construction defines a map $e:\mathcal{C}_X\rightarrow \mathcal{M}_X$. 
Under the hypothesis of the Proposition, this map is smooth, in fact its fiber over $E\in \mathcal{M}_X$ is the open subset of the projective space $\P(H^0(E))$ corresponding to sections $s$ vanishing  along a curve of $\mathcal{C}_X$.

\medskip	
\emph{Remark} $3.-$ Under the hypothesis of Proposition \ref{deform}, we have $H^1(\mathcal{I}_CL)=0$, hence, in view of (\ref{ext}),  $H^1(E)=0$. 
Assume moreover $H^3(E^*)=0$ (or, equivalently, $H^0(E\otimes K_X)=0$). 
From the exact sequence (\ref{end}) we get an isomorphism 
$H^2(\mathcal{E}nd(E))\iso  H^2(\mathcal{I}_CE)$. From the exact sequence $0\rightarrow \mathcal{I}_CE\rightarrow E\rightarrow E_{|C}\rightarrow 0$ we get an injection of $H^1(E_{|C})=H^1(N_C)$ into $H^2(\mathcal{I}_CE)$. Hence the vanishing of $H^1(N_C)$ is in fact \emph{equivalent} in this case to that of  $H^2(\mathcal{E}nd(E))$.

\medskip	
\section{The examples: set-up}
The rest of the paper is devoted to examples. We will not hesitate to make strong hypotheses to simplify the exposition -- we encourage the reader to explore more general situations.

Thus we will only consider Fano threefolds whose Picard group is cyclic and generated by a very ample line bundle $\O_X(1)$.  Then $K_X=\O_X(-i)$, where $i$ is the \emph{index} of $X$.
Leaving aside the projective space and the quadric, we get two series:

$\bullet$ Index 2: $X_d\subset \P^{d+1}$, for $3\leq d\leq 5$;

$\bullet$ Index 1: $X_{2g-2}\subset \P^{g+1}$, with $3\leq g\leq 10$ or $g=12$.

We will choose our curve $C$ so that $L=\O_X(1)$ or $\O_X(2)$, and $K_C=\O_C(j)$ with $j=-1,0$ or $1$; moreover $C$ will be projectively normal. This implies that the surjectivity assumptions in Propositions (\ref{main}) and (\ref{deform}) are automatically satisfied. 

We will assume  $S\not\supset C$, and that $\Pic(S)$ \emph{is generated by $\O_S(1)$}. This is the case if $S$ is sufficiently general \cite[Theorem 3.33]{V}; it will allow us to deal only with \emph{stable} vector bundles, thanks to the following lemma:
\begin{lem}\label{stable}
If $L=\O_X(1)$, or $L=\O_X(2)$ and $C$ is not contained in a hyperplane, $E$ and $E_S$ are stable.
\end{lem}
\pr Since $\Pic(X)=\Z\cdot \O_X(1)$,   $E$  is stable if and only if $H^0(E(-1))=0$; if $L=\O_X(1)$ this is automatic, if $L=\O_X(2)$ this is equivalent to $H^0(\mathcal{I}_C(1))=0$, which means that $C$ is not contained in a hyperplane. 

Put $L=\O_X(j)$. Let $Z$ be the finite subscheme $C\cap S$. Restricting (\ref{ext}) to $S$ gives an exact sequence
\[0\rightarrow \O_S \rightarrow E_S \rightarrow \mathcal{I}_Z(j)\rightarrow 0\,.\]
If $j=1$, $E_S$ is stable; if $j=2$, $E_S$ is stable if $H^0(\mathcal{I}_Z(1))=0$. Using the exact sequence \allowbreak$0\rightarrow \mathcal{I}_C(-i)\rightarrow \mathcal{I}_C\rightarrow \mathcal{I}_{Z}\rightarrow 0$, we see that if   $C$ is not contained in a hyperplane, it suffices to show $H^1(\mathcal{I}_C(1-i))=0$. But this follows from the exact sequence $0\rightarrow \mathcal{I}_C(1-i)\rightarrow \O_X(1-i)\rightarrow \O_C(1-i)\rightarrow 0$.\qed

\smallskip	
We will always assume that the hypotheses of the Lemma are satisfied, and 
 we will take for $\mathcal{M}_S$ the moduli space of stable rank 2 vector bundles with $c_1=[L]$ and $c_2=\deg(C)$. 

\medskip	
\section{Rational curves}
We first consider the case where the curve $C$ is rational.
 The condition $(K_X\otimes L)_{|C}=K_C$ with $L$ ample imposes 
 $\deg(C)=2$, $L=\O_X(1)$, $K_X=\O_X(-2)$. Thus $C$ is a conic,  $X$ has index 2, $S$ is the intersection of $X$ with a quadric $Q$. We have an
 extension $0\rightarrow \O_X\rightarrow E\rightarrow \mathcal{I}_C(1)\rightarrow 0$. 
 
Let $h\in H^2(X,\Z)$ be the class of a hyperplane section. We have $c_1 (E)=h$ and $c_2(E)\cdot h= C\cdot h=2$, hence $\Delta _E\cdot h=8-d$.
 Assuming $H^1(N_C)=0$, we deduce from (\ref{RR}) $\dim \mathcal{M}_X=5-d$.
   
 \medskip	  
\noindent$\boldsymbol{d=3}$. We let $X\subset\P^4$ be a cubic threefold, and $S=X\cap Q$ for a quadric $Q$.
 \begin{prop}
$1)$ The conics $C\subset X$ satisfy $H^1(N_C)=0$.

$2)$ The moduli space $\mathcal{M}_X$ is isomorphic to the Fano surface of lines contained in $X$.

$3)$  The moduli space $\mathcal{M}_S$ is isomorphic to the Hilbert square $S^{[2]}$ of $S$.

$4)$ With the above identifications, the map $\operatorname{res}: \mathcal{M}_X\rightarrow \mathcal{M}_S$ associates to a line $\ell$ the length 2 subscheme $\ell\cap Q$ of $S$. It is an isomorphism onto a closed Lagrangian submanifold of $S^{[2]}$.
\end{prop}

\pr 1) Let $\mathcal{C}_X$ be the Hilbert scheme of (possibly degenerate) conics contained in $X$, and let $F_X$ be the Fano surface of lines. Associating  to  a conic its residual line gives a morphism $\rho :\mathcal{C}_X\rightarrow F_X$, which realizes $\mathcal{C}_X$ as a $\P^2$-bundle over $F_X$. In particular, $\mathcal{C}_X$ is smooth, of dimension 4. Therefore $h^0(N_C)=4$; by Riemann-Roch this implies $H^1(N_C)=0$. 

\smallskip	
$2)$ Let $E\in\mathcal{M}_X$. 
We have $h^0(E)=3$; each nonzero section defines a  conic in $X$. We claim that all these conics have the same residual line $\ell$. The most economical way to prove this is to use the relation $c_2(E)=[C]=h^2-\ell$ in the Chow group $CH^2(X)$. If $C'$ is another conic defined by a section of $E$, with residual line $\ell'$, we have $\ell'-\ell=0$ in $CH^2(X)$, hence in the intermediate Jacobian $JX$; but the Abel-Jacobi map $\ell'\mapsto \ell'-\ell$ embeds $F_X$ into $JX$, hence $\ell'=\ell$. 

It follows that the map $e:\mathcal{C}_X\rightarrow \mathcal{M}_X$ (\S 2) factors as $\mathcal{C}_{X}\rightarrow F_X\xrightarrow{\ u\ }\mathcal{M}_X$, where $u$ is bijective, hence an isomorphism since $F_X$ and $\mathcal{M}_X$ are smooth.

\smallskip	
$3)$ Consider the moduli space of stable rank 2 bundles $F$ on $S$ with $c_1=h$, $c_2=4$. Riemann-Roch gives $h^0(F)\geq 3$; since $\Pic(S)=\Z$, a  nonzero section $s$ of $F$ vanishes along a finite subscheme $Z$ of length 4. Thus we have an exact sequence 
\[0\rightarrow \O_S \xrightarrow{\ s\ } F \rightarrow \mathcal{I}_Z(1)\rightarrow 0 \, ,\]with $h^0(\mathcal{I}_Z(1))\geq 2$ (and actually $=2$, since $Z$ cannot be contained in a line), so that $Z$ is contained in a 2-plane. 
This plane meets $S$ along a finite subscheme of length 6, contained in the conic $Q\cap\lr{Z}$; let $\mathfrak{z}$ be the residual subscheme of $Z$ in that conic. As before, we claim that  $\mathfrak{z}$ does not depend on the choice of the section $s$.

The exterior product $\varphi :\bigwedge^2H^0(F)\rightarrow H^0(\det F)$ is injective: indeed, since $h^0(F)=3$, any element of $\bigwedge^2H^0(F)$ is of the form $s\wedge t$, with $s,t\in H^0(F)$; if $\varphi (s\wedge t)=0$,  the exact sequence $0\rightarrow \O_S\xrightarrow{\ s\ } F \xrightarrow{\ \wedge s\ } \mathcal{I}_Z(1)\rightarrow 0$ shows that $t\in \C s$. 

Thus $\im \varphi $ has codimension $2$ in $H^0(\O_S(1))$, hence consists of the sections vanishing along a line $\ell\subset \P^4$. Therefore the evaluation homomorphism $H^0(F)\otimes_{\C}\O_S \rightarrow F$ is surjective outside the finite set $\ell\cap S$, with kernel $(\det F)^{-1}=\O_S(-1)$. Dualizing, we get an exact sequence:
\[0\rightarrow F^*\rightarrow  H^0(F)^*\otimes_{\C}\O_S\rightarrow \mathcal{I}_{A}(1)\rightarrow 0\, ,\]where $A$ is a subscheme  supported in $\ell\cap S$. Computing $c_2$ gives that $A$ has length 2. 

When $t$ runs through $H^0(F)$, the hyperplanes  $\varphi (s\wedge t)=0$ form a pencil, whose intersection is $\lr{Z}$; thus $\lr{Z}$ contains $\ell$ and therefore $A$.  A general section of $F$ does not vanish at a point $a\in A$, since otherwise all the hyperplanes containing $\ell$ would intersect $S$ with multiplicity $\geq 2$ at $a$. Therefore $A$ is indeed the residual subscheme of all zero loci of sections of $F$. 

\smallskip	
The zero locus $Z$ of a nonzero section of $s$ is   a local complete intersection, contained in a 2-plane, and
 has the \emph{Cayley-Bacharach property}  \cite[Theorem 5.1.1]{Hu-L}: no length 3 subscheme $Z'\subset Z$ is contained in a line. Let $\mathcal{H}$ be the locally closed subscheme of $S^{[4]}$ parametrizing subschemes with these properties, and let $Z\in \mathcal{H}$. Serre duality  provides an isomorphism $\Ext^1_S(\mathcal{I}_Z(1),\O_S) \iso H^1(\mathcal{I}_S(1))^*$; using the exact sequence $0\rightarrow \mathcal{I}_Z(1)\rightarrow \O_S(1)\rightarrow \O_Z(1)\rightarrow 0$ and $h^0(\mathcal{I}_Z(1))=2$, we get that this space has dimension 1. Therefore 
 there exists a unique nontrivial extension $0\rightarrow \O_S\rightarrow F\rightarrow \mathcal{I}_Z(1)\rightarrow 0$, and the Cayley-Bacharach property ensures that $F$ is locally free (\emph{loc. cit.}). We get in this way a morphism $\mathcal{H}\rightarrow \mathcal{M}_S$; from the above it factors as $\mathcal{H}\rightarrow S^{[2]}\rightarrow \mathcal{M}_S$, where the first map associates to a subscheme $Z$ its residual subscheme in the conic $\lr{Z}\cap S$. The second map is bijective, hence an isomorphism.

\smallskip	
$4)$ Let $C$ be a conic in $X$, and let $E$ be the associated vector bundle. The restriction $E_S$ is associated to the subscheme $Z:=C\cap Q$ of $S$;  the residual subscheme of $Z$ in $\lr{Z}\cap Q=\lr{C}\cap Q$ is $\ell\cap Q$, where $\ell$ is the residual subscheme of $C$ in $\lr{C}\cap X$. It is clear that the map $\ell\mapsto \ell\cap Q$ is an embedding; by $1)$ and Tyurin's theorem, its image is a Lagrangian submanifold of $S^{[2]}$.\qed

\medskip	
\noindent$\boldsymbol{d=4}$.	
We consider now the degree 4 case. Then $X$ (resp.\ $S$) is a complete intersection of $2$ (resp.\ $3$) quadrics in $\P^5$. The quadrics in $\P^5$ containing $S$ form a net $\Pi\cong\P^2$, and those containing $X$ form a line $\ell\subset \Pi$. Let $\Delta \subset\Pi$ be the degree 6 discriminant curve  parametrizing singular quadrics. Since $\Pic(S)=\Z$, 
 $\Delta $ is smooth: this is equivalent to say that every quadric in $\Delta $ has rank $5$ \cite[Proposition 1.2]{BP}, and this holds because otherwise $S$ would be contained in a quadric of rank $\leq 4$, and the 3-planes of that quadric would cut a degree 4 curve on $S$. The double covering $\pi :\hat{S}\rightarrow \Pi$   branched along $\Delta $ parametrizes the pairs $(Q,\sigma )$, where $Q\in \Pi$ and $\sigma $ is one of the two families of 2-planes contained in $Q$. The surface $\hat{S}$ is a K3 surface, and $\Gamma :=\pi ^{-1}(\ell)$ is a curve of genus 2, which is well known to play a fundamental role in the geometry of $X$ (for instance, the intermediate Jacobian of $X$ is the Jacobian of $\Gamma $).

\begin{prop}\label{d4}
$1)$ The conics $C\subset X$ satisfy $H^1(N_C)=0$. 

$2)$ The moduli space $\mathcal{M}_X$ is isomorphic to $\Gamma $.

$3)$ The moduli space $\mathcal{M}_S$ is isomorphic to $\hat{S}$, and the restriction map $\mathcal{M}_X\rightarrow \mathcal{M}_S$ corresponds to the embedding $\Gamma \hookrightarrow \hat{S}$.
\end{prop}

\pr Let $C$ be a conic contained in $X$, and let $\lr{C}$ be the 2-plane spanned by $C$. 
 Since every quadric of $\ell$ contains $C$,  the plane spanned by $C$ must be contained in one (and only one) quadric $Q$ of $\ell$; we get a map $\mathcal{C}_X\rightarrow \Gamma $ by associating to $C$ the quadric $Q$ and the family $\sigma $ of 2-planes in $Q$ containing $\lr{C}$. 

The vector bundle $E$ associated to $C$ admits a simple description. If the quadric $Q$ is smooth, we identify it to the Grassmannian $\G(2,4)$, in such a way that the zero loci of the nonzero sections of the universal quotient bundle $G$ are the 2-planes of $\sigma $;   then $E=G_{|X}$.

If $Q$ is singular, it is a cone with vertex $v$ over a smooth quadric $Q_0\subset \P^4$, which we  view as a hyperplane section of $\G(2,4)$. Let $G'$ be the pull back to $Q\smallsetminus\{v\} $ of $G_{|Q_0}$; the zero loci of the nonzero sections of  $G'$ are the 2-planes of $\sigma $ minus $v$. 
 Therefore $E=G'_{|X}$. In each case   the zero loci of the nonzero sections of $E$ are the conics $P\cap X$ for $P\in\sigma $. 

In each case, the vector bundle $E$ is globally generated, and therefore $N_C=E_{|C}$ is globally generated. Since $H^1(C,\O_C)=0$, this implies $H^1(N_C)=0$. Therefore $\mathcal{M}_X$ is smooth; the map $\mathcal{C}_X\rightarrow \Gamma $ factors as $\mathcal{C}_{X}\xrightarrow{\ e\ } \mathcal{M}_X \xrightarrow{\ u\ }\Gamma $, where $u$ is bijective, hence an isomorphism. 

Finally let us consider the moduli space $\mathcal{M}_S$ of stable bundles on $S$ with $c_1=h$, $c_2=4$. By \cite[Example 0.9]{M1}, it is isomorphic to $\hat{S}$: as above, if $Q$ is a smooth quadric in $\Pi$, we identify $Q$ to $\G(2,4)$
so that each nonzero section of $G$ vanish along a 2-plane of $\sigma $, and take $F=G_{|S}$. By construction the restriction map $\mathcal{M}_X\rightarrow \mathcal{M}_S$ coincides with the natural embedding $\Gamma \hookrightarrow \hat{S}$.\qed

\medskip	
Note that any curve in a K3 surface is Lagrangian, so Tyurin's theorem does not give any information in this case. 

\medskip	
\noindent$\boldsymbol{d=5}$. Finally we look at threefolds $X\subset \P^6$ of degree 5 and index 2. Recall that such a threefold is the section by a $6$-plane of the Grassmannian $\G(2,5)\subset\P^9$.
 
 \begin{prop}
The moduli spaces $\mathcal{M}_X$ and $\mathcal{M}_S$ consist of one reduced point, which is the restriction of the universal quotient bundle on $\G(2,5)$.
\end{prop}
\pr A conic  $C\subset \G(2,5)$ corresponds to a  surface of degree 2 in $\P^4$, which is necessarily contained in a $3$-plane of $\P^4$; therefore $C$ is contained in a sub-Grassmannian $\G(2,4)\subset \G(2,5)$. This $\G(2,4)$ is the zero locus of a section $s$ of the universal quotient bundle $G$ on $\G(2,5)$; the restriction of $s$ to  $X$ vanishes along the intersection of $\G(2,4)$ with a 2-plane. Since $X$ does not contain a plane, the zero locus of $s_{|X}$ is $C$. Therefore $E=G_{|X}$. Again since $G$ is globally generated, $N_C$  is globally generated, thus $H^1(N_C)=0$, and $\mathcal{M}_X$ consists of the reduced point $\{G_{|X}\} $.

Since $\mathcal{M}_S$ is zero-dimensional, it also consists of a unique reduced point \cite[\S 3]{M2}, which is given by $G_{|S}$.\qed

\medskip	
\section{Elliptic curves, index 2}
After rational curves, the next case is elliptic curves; we must take $L=\O_X(i)$, where $i$ is the inde of $X$. Let us first look at the case where $X$ has index 2,   that is, $X_d\subset\P^{d+1}$ with $d=3,4$ or $5$. The surface $S$ is the intersection of $X$ with a quadric $Q$.

The vector bundle $E$ associated to $C$ is  stable if and only if  $C$ spans $\P^{d+1}$ (Lemma \ref{stable});
 it is therefore natural to look at  \emph{normal elliptic curves} $C_{d+2}\subset \P^{d+1}$, embedded by a complete linear series of degree $d+2$.
\begin{prop}
$1)$ Every $X_d$ contains a normal elliptic curve $C_{d+2}\,\cdot $

$2)$ The associated rank 2 vector bundle $E$ is a \emph{Ulrich bundle} -- that is, $H^{\pp}(E(-i))=0$ for $i=1,2,3$.

$3)$  We have $\dim\mathcal{M}_X=5$ and $\dim\mathcal{M}_S=10$; $\mathcal{M}_S$ is birational to the O'Grady manifold $\operatorname{OG} _{10}$. The restriction map $\mathcal{M}_X \rightarrow \mathcal{M}_S$ is injective, and a local isomorphism onto a lagrangian subvariety. 
\end{prop}
\pr All these assertions are proved in \cite[\S6]{U}, except the injectivity of $\operatorname{res} $. Let $E,F$ be two elements of $\mathcal{M}_X$; since $E$ is a Ulrich bundle, we have a presentation
$\O_{\P}^b(-1)\rightarrow \O_{\P}^a \rightarrow E\rightarrow 0$, and, by restriction to $Q$, $\O_{Q}^b(-1)\rightarrow \O_{Q}^a \rightarrow E_S\rightarrow 0$. Applying the functor $\Hom(-,F)$ we get a commutative diagram
\[\xymatrix{ 0\ar[r]& \Hom(E,F)\ar[r]\ar[d]^{\alpha } & H^0(F)^a \ar[r]\ar[d]^{\beta } & H^0(F(1))^b\ar[d]^{\gamma } \\
0\ar[r]& \Hom(E_{S},F_{S})\ar[r]  & H^0(F_{S})^a \ar[r]  & H^0(F_{S}(1))^b}\]
Using the exact sequence $0\rightarrow F(-2)\rightarrow F\rightarrow F_S\rightarrow 0$ and the vanishing of $H^{\pp}(F(-1))$ and $H^{\pp}(F(-2))$ we see that $\beta $ and $\gamma $ are bijective, hence $\alpha $ is bijective. Thus if $E_S$ and $F_S$ are isomorphic, there is a nonzero homomorphism from $E$ to $F$, which must be an isomorphism since $E$ and $F$ are stable.\qed

\medskip	
\noindent\textbf{Example :} $\boldsymbol{d=3}$

  This case  has two interesting features, which are treated in detail in \cite{Bcub}:

a) The moduli space $\mathcal{M}_X$ is birational to the intermediate Jacobian $JX$. 

b) Let us fix the K3 surface $S$ (a $(2,3)$ complete intersection in $\P^4$). The projective space $\Pi$ of cubic hypersurfaces containing $S$ has dimension $5$; there is a rational Lagrangian fibration $h:\mathcal{M}_S\dasharrow \Pi$, whose fiber at a general cubic $X$ is isomorphic to $\mathcal{M}_X$ -- hence birational to $JX$. 

Suppose $S$ is given by $Q=F=0$, with $\deg(Q)=2$, $\deg(F)=3$. The elements of $\Pi$ are the cubics $aF+LQ=0$, with $a\in\C$ and $L\in H^0(\P^4,\O_{\P^4}(1))$; for $a\neq 0$, this cubic 
 can be identified with the  section by the hyperplane $L=aT$ of the nodal cubic \emph{fourfold} $V\subset \P^5$ defined by $F+TQ=0$    (here $T$ is a new coordinate added to the coordinates on $\P^4$).  In other words, we can view $\Pi$ as the dual of $\P^5$; then the fiber of the map $h:\mathcal{M}_S\dasharrow (\P^5)^*$ at a general hyperplane $H$ is birational to the intermediate Jacobian $J(V\cap H)$. This should be compared with \cite{LSV}, where the authors construct a Lagrangian  fibration $h:\operatorname{OG}_{10}\rightarrow (\P^5)^* $ for a very general cubic fourfold $V\subset \P^5$, such that $h^{-1}(H)=J(V\cap H)$ for a general $H\in(\P^5)^*$. 
 
In our case, we do not know whether there exists a projective holomorphic symplectic manifold $\mathcal{M}'$ and a birational map $\mathcal{M}'\bir \mathcal{M}_S$ such that the composition $\mathcal{M}'\bir \mathcal{M}_S\dasharrow (\P^5)^*$ is everywhere defined.


\section{Elliptic curves, index 1}

We can also perform the Serre construction from elliptic curves lying on an index 1 Fano threefold -- this time we do not need to assume that $C$ spans the projective space. We will work out one case, which gives an example where the restriction map is \emph{not} injective.

We take for  $X$  a complete intersection of 3 quadrics in $\P^6$, and for $C\subset X$ a normal elliptic curve in $\P^3$,   complete intersection of 2 quadrics. We have an extension \[0\rightarrow \O_X\rightarrow E\rightarrow \mathcal{I}_C(1)\rightarrow 0\,.\]

Recall that the quadrics containing $X$ form a 2-dimensional projective space $\Pi=\abs{\mathcal{I}_X(2)}$, and that the discriminant curve $\Delta \subset\Pi$ parametrizing singular quadrics is a degree 7 nodal curve. We will assume for simplicity that $\Delta $ \emph{is smooth} -- this is equivalent to say that all quadrics in $\Delta $ have rank $6$, and this holds when $X$ is general \cite[Proposition 1.2]{BP}. Then the two families of $3$-planes contained in any quadric of $\Delta $ define an \'etale double covering $\rho  :\tilde{\Delta }\rightarrow \Delta  $.

The K3 surface $S$ is the intersection of $X$ with a hyperplane $H\subset\P^5$. The quadrics of $H$ containing $S$ are again parametrized by $\Pi$, via the restriction map $Q\mapsto Q\cap H$. The discriminant curve $\Delta _H$ has now degree $6$. As in \S 4 (case $d=4$), we consider the double covering $\pi :\hat{S}\rightarrow \Pi$ branched along $\Delta _H$. We define a map $r:\tilde{\Delta }\rightarrow \hat{S} $ by associating to a pair $(Q,\sigma )$ in $\tilde{\Delta } $ the pair $(Q\cap H,\sigma _H)$, where $\sigma _H$ is the family of 2-planes $P\cap H$ for $P\in \sigma $.

\begin{prop}
$1)$ The quartic elliptic curves $C\subset X$ satisfy $H^1(N_C)=0$. 

$2)$ The moduli space $\mathcal{M}_X$ is isomorphic to $\tilde{\Delta } $.

$3)$ The moduli space $\mathcal{M}_S$ is isomorphic to $\hat{S}$.

$4)$ The restriction map $\operatorname{res}: \mathcal{M}_X\rightarrow \mathcal{M}_S$ is identified with $r:\tilde{\Delta} \rightarrow \hat{S} $. This map factors as   \break  $\tilde{\Delta }\xrightarrow{\ n\ } \pi ^{-1}(\Delta ) \hookrightarrow \hat{S}$; the curve $\pi ^{-1}(\Delta )$ is singular along $\pi ^{-1}(\Delta \cap \Delta _H)$, and $n$ is its normalization. In particular, $\operatorname{res} $ is not injective. 
\end{prop}

\pr Let $P$ be the 3-plane spanned by $C$.	
The kernel of
the restriction map $H^0(\P^6,\mathcal{I}_X(2))\rightarrow \allowbreak H^0(P,\mathcal{I}_C(2))$ has dimension $\geq 1$, and actually $1$ since otherwise $X$ would contain a quadric surface. Thus there is a unique quadric $Q_1\in\Pi$  which contains $P$. 

1) We can suppose that $C$ is defined by $X=Y=Z=R=S=0$, where $R$ and $S$ are quadratic forms. 
Then $X$ is defined by 3 equations $Q_i=0$, with
$Q_i=XL_i+YM_i+ZN_i+ a_iR+b_iS$,   $L_i\in H^0(\P^6,\O_{\P}(1))$ and $a_i,b_i\in\C$ $(i=1,2,3)$, $a_1=b_1=0$.

The exact sequence of normal bundles for $C\subset X\subset\P^6$ reads
\[0\rightarrow N_C \rightarrow \O_C(1)^3\oplus \O_C(2)^2 \xrightarrow{\ M\ } \O_C(2)^3\rightarrow 0 \, ,\]
where $M$ is the $3\times 5$ matrix with rows $(L_i,M_i,N_i,a_i,b_i)$.

Suppose $H^1(N_C)\neq 0$. Then by Serre duality 
 $H^0(N_C^*)\cong H^0(N_C(-1))$ is nonzero, so
  there exist scalars $\ell,m,n$, not all zero, such that 
$\ell L_1+mM_1+nN_1=0$ in $H^0(C,\O_C(1))$; that is, the linear forms $X,Y,Z,L_1,M_1,N_1$ are linearly dependent. But this implies that the quadric $Q_1$ has rank $\leq 5$, contradicting the hypothesis. 

\medskip	
2) The quadric $Q_1$ is a cone over a smooth quadric $G\subset\P^5$, with vertex $v$; the 3-plane $P$ is spanned by $v$ and a 2-plane $P_0\subset G$. As in the proof of Proposition \ref{d4}, 
 identifying $G$ with the Grassmannian of lines in $\P^3$ gives  
 a  rank 2 vector bundle $E_0$ on $G$, with a section $s_0$ vanishing along $P_0$. 
 Pulling back $(E_0,s_0)$ to $Q_1\smallsetminus\{v\} $ and restricting to $X$ gives a rank 2 bundle  on $X$, with a section  vanishing along $C$: this is our bundle $E$. Varying $s_0$ gives all 3-planes contained in $Q_1$, so the map $\mathcal{C}_X\rightarrow \tilde{\Delta } $ factors through an isomorphism $\mathcal{M}_X\rightarrow \tilde{\Delta } $.
 
 \medskip	
 3) We have already seen that the moduli space $\mathcal{M}_S$ is isomorphic to $\hat{S}$ (Proposition \ref{d4}); 
 by construction the restriction map $\operatorname{res}:\mathcal{M}_X\rightarrow \mathcal{M}_S $ corresponds to the map $r:\tilde{\Delta }\rightarrow \hat{S} $ given by $(Q,\sigma )\mapsto (Q\cap H, \sigma _H)$. 
 
 \smallskip	
\raisebox{-5pt}{4) Consider the commutative diagram}  \qquad\xymatrix@M=5pt{\tilde{\Delta }\ar[r]^n\ar[dr] _{\rho } & \pi ^{-1}(\Delta )\ar@{^{(}->}[r]\ar[d]^{\pi }& \hat{S}\ar[d]^{\pi }\\
& \Delta\ar@{^{(}->}[r] & \Pi}

\smallskip	
Put $\Sigma :=\Delta \cap\Delta _H$. Above $\Delta  \smallsetminus \Sigma $, $\rho $ and $\pi $ are \'etale double coverings, hence $n$ is an isomorphism. At a point $Q$ of $\Sigma $, $n$ maps the two points of $\rho ^{-1}(Q)$  to the unique point of $\pi ^{-1}(Q)$.  Therefore 
$\pi ^{-1}(\Delta )$ is singular along $\pi ^{-1}(\Sigma )$,
$n$ is its normalization,    and $\operatorname{res} $ is not injective.\qed

\medskip	
\section{Canonical curves}

Another way of making sure that $K_{C}$ is the restriction of a line bundle on $X$ is to take for $C$ a canonical curve, and $X=X_g\subset \P^{g+1}$ of index 1. 
Then $L=\O_X(2)$, so we have an exact sequence
\[0\rightarrow \O_{X}\rightarrow E\rightarrow \mathcal{I}_C(2)\rightarrow 0\,.\]In view of  Lemma \ref{stable}, we will take for $C\subset \P^{g+1}$ a curve of genus $g+2$ embedded by its canonical system. 

\begin{prop}
Let $X\subset \P^{g+1}$ be a  Fano threefold of index $1$. If $g=3$ we assume that $X$ is general. 
 
$1)$ $X$ contains a canonical curve $C$ of genus $g+2$, satisfying
 $H^1(N_C)=0$.

$2)$ The moduli space $\mathcal{M}_S$ is birational to $\operatorname{OG}_{10} $; the restriction map $\operatorname{res}: \mathcal{M}_X\rightarrow \mathcal{M}_S $ is a local isomorphism into a Lagrangian subvariety of $\mathcal{M}_S$.
\end{prop}

\pr We first consider the case $g\geq 4$. We will rely heavily on a result of \cite{B-F}: $X$ carries a stable rank 2 vector bundle $F$ with $\det(F)=\O_X$ and $c_2(F)\cdot h=4$, which satisfies $H^2(\mathcal{E}nd(F))=0$.  
The vector bundle $F$ is constructed as a flat deformation of a torsion free coherent sheaf $\mathcal{F}$ which fits into an extension
\[0\rightarrow \mathcal{I}_A\rightarrow \mathcal{F}\rightarrow \mathcal{I}_B\rightarrow 0\,,\]where $A$ and $B$ are two general conics in $X$. Since $g\geq 4$ the 2-plane spanned by $A$ (or $B$) does not contain any other point of $X$; therefore $\mathcal{I}_A(1)$ and $\mathcal{I}_B(1)$ are globally generated. Since $H^1(\mathcal{I}_A(1))=0$, it follows that $\mathcal{F}(1)$ is globally generated. We have $h^{i}(\mathcal{F}(1))=0$ for  $i\geq 1$, hence $h^0(F(1))=h^0(\mathcal{F}(1))$ for $F$ general enough; therefore $E:=F(1)$ is globally generated.

The zero locus $C$ of a general section of $E$  is a smooth curve, with normal bundle $N_{C}=E_{|C}$. The adjunction formula gives $K_C=\O_C(1)$; the exact sequence $0\rightarrow \O_X(-1)\rightarrow F\rightarrow \mathcal{I}_C(1)\rightarrow 0$ gives $h^0(\mathcal{I}_C(1))=h^1(\mathcal{I}_C(1))=0$, hence the restriction map from $H^0(\P^{g+1},\O_{\P})(1)= H^0(X,\O_X(1))$ to $H^0(C,K_C)$ is an isomorphism. Therefore $C$ is a curve of genus $g+2$, canonically embedded in $\P^{g+1}$. We have $H^0(E\otimes K_X)=H^0(F)=0$, hence $H^1(N_C)=0$ by  Remark 3. This proves $1)$ in the case $g\geq 4$. Then Tyurin's theorem applies, showing that $\operatorname{res}: \mathcal{M}_X\rightarrow \mathcal{M}_S$ is a local isomorphism into a Lagrangian subvariety of $\mathcal{M}_S$; since 
the vector bundle $E_S(-1)=F_{|S}$ has $c_1=0$, $c_2=4$, $\mathcal{M}_S$ is birational to $\operatorname{OG}_{10} $.

\medskip	
We now consider the case $g=3$: we want to prove that a general quartic hypersurface in $\P^4$ contains an intersection of 3 quadrics. 
Consider the map $p:H^0(\P^4, \O_{\P}(2))^6\rightarrow H^0(\P^4, \O_{\P}(4))$ given by $p(q_0,\ldots ,q_5)=\sum\limits_{i=0}^2 q_iq_{i+3}$. Its differential at $(q_i)$ associates to $(r_i)\in H^0(\P^4, \O_{\P}(2))^6$ the   form $\sum q_ir_i$ (after permuting the indices). 
This differential  is surjective  for $(q_i)$ general: it suffices to prove it for one particular sextuple $(q_i)$;  taking $q_i=X_i^2$ for $0\leq i\leq 4$ and $q_5=X_0X_1+X_2X_3$, one sees easily that every degree 4 monomial belongs to the image.  Thus $p$ is dominant, so that a general quartic in $\P^4$ admits an equation of the form $\sum\limits_{i=0}^2 q_iq_{i+3}=0$, hence contains the canonical curve $C$ of genus 5 defined by $q_0=q_1=q_2=0$.

The exact sequence of normal bundles for $C\subset X\subset \P^4$ becomes here
\[0\rightarrow N_C \rightarrow \O_C(2)^3 \xrightarrow{\ (q_3,q_4,q_5) \ } \O_C(4)\rightarrow 0\,.\]Suppose $H^1(N_C)\neq 0$.
By Serre duality, $H^0(N_C^*(1))=H^0(N_C(-1))\neq 0$. In view of 
the above exact sequence, this means that there exist a nonzero triple of 
 linear forms $(\ell_3,\ell_4,\ell_5)$ on $\P^4$ such that $\sum\limits_{j=3}^5 \ell_jq_j=0 $ in $H^0(C, \O_C(3))$; in other words, a nonzero sextuple $(\ell_0,\ldots ,\ell_5)$ such that $\sum\limits_{i=0}^5 \ell_iq_i=0$. But for a general choice of the $q_i$ such a sextuple does not exist: again taking $q_i=X_i^2$ for $0\leq i\leq 4$ and $q_5=X_0X_1+X_2X_3$ does the job. This proves $1)$ in this case, and $2)$ follows as above.\qed

\bigskip	
\noindent\textbf{Example: $\boldsymbol{g=4}$} 

When $g=4$, $X$  is a $(2,3)$-complete intersection in $\P^5$. In this case there is a simple way to find a genus 6 canonical curve in $X$: in the space of cubic fourfolds,
those which contain a Del Pezzo surface $S_5\subset \P^5$ of degree 5 form a hypersurface  ($\mathcal{C}_{14}$ in the notation of \cite{H}). Since the space of cubics containing $X$ is 6-dimensional, a general $X$ can be written $X=V\cap Q\subset \P^5$, where $V$ is a cubic containing a Del Pezzo quintic  $S_5$ and $Q$ is a quadric. Then $X$ contains the curve $S_5\cap Q$, which is canonical of genus 6. 
We do not know whether the corresponding component   $\mathcal{M}_X$ of the moduli space obtained in this way is the same as the component described in \cite{B-F}. But we still have the required vanishing:

\begin{lem}
For $X$ general and all canonical curves $C\subset X$, we have $H^1(N_C)=0$.
\end{lem}
\pr We will denote by $H_c$ be the Hilbert scheme of canonically embedded curves of genus $6$    $C\subset\P^{5}$, by $H_x$ the Hilbert scheme of  smooth $(2,3)$-complete intersections  in $\P^5$,   and by $H_{c,x}$ the nested Hilbert scheme of pairs $C\subset X\subset \P^{g+1}$. Consider the diagram:
\[\xymatrix{ & H_{c,x} \ar[dl]_p\ar[dr]^q& \\ H_c && H_x\ .
}\]We first observe that our three Hilbert schemes are smooth. This is clear for $H_{x}$, and well-known for $H_{c}$ -- in fact we have $H^1(C,N_{C/\P})=0$. 
Finally, 
 recall that the tangent space to $H_{c,x}$ at $(C,X)$ fits into a cartesian diagram
\[\xymatrix{ T_{(C,X)} H_{c,x}\ar[r]^<<<<<<<<<{T(q)}\ar[d]_{T(p)} & T_X(H_x)=H^0(X,N_{X/\P})\ar[d]^r\\
T_{C}(H_c)=H^0(C, N_{C/\P})\ar[r]^>>>>>>>{u } &  H^0(C, N_{X/\P}{}^{}_{|C})}
\](see for instance \cite[Lemma 8.8]{ACG}). Since  $N_{X/\P}=\O_X(2)\oplus \O_X(3)$ and $C$ is projectively normal,  the restriction homomorphism $r$ is surjective, and therefore $T(p)$ is surjective; hence $p$ is smooth, and $H_{c,x}$ is smooth. 
Therefore  $q$ is smooth along a general fiber $q^{-1}(X)$,  hence $T(q)$ is surjective at all points  $C\subset X$; since $r$ is surjective, this implies that $u$ is surjective. 
From the exact sequence of normal bundles for $C\subset X\subset \P$ we get an exact sequence
$H^0(C, N_{C/\P}) \xrightarrow{\ u \ } H^0(C, N_{X/\P}{}^{}_{|C})\rightarrow \allowbreak H^1(N_C) \rightarrow 0$, hence $H^1(N_C)=0$. \qed

\medskip	

Thus again the restriction map $\mathcal{M}_X\rightarrow \mathcal{M}_S$ is a Lagrangian immersion.
We can say a little more in this case. We have realized the canonical curve $C$ as $S_5\cap Q$, where $S_5$ is a quintic Del Pezzo surface contained in $V$. Now the Serre construction determines uniquely a rank 2 vector bundle $F$ on $V$ 
and an extension $0\rightarrow \O_V\rightarrow F\rightarrow \mathcal{I}_{S_5}(2)\rightarrow 0$; 
the restriction of $F$ to $X$ is the bundle $E$ associated to $C$. 

Recall from \cite{Bdet} that $F$ fits into an exact sequence \begin{equation}\label{M}
0\rightarrow \O_{\P^5}(-1)^6\xrightarrow{\ M\ }\O_{\P^5}^6\rightarrow F\rightarrow 0\, ,
\end{equation}
 where $M$ is a skew-symmetric $6\times 6$ matrix of linear forms; this implies that the cubic form defining $V$ is the Pfaffian of the  skew-symmetric  matrix $M$. For a general
pfaffian cubic $V$ the vector bundle $F$ is uniquely determined \cite[Proposition 9.2 (b)]{Bdet}. Let $\mathfrak{P}\,(\cong\P^6)$ be the space of cubics in $\P^5$ containing $X$, and
let $\mathpzc{Pf}$ be the hypersurface of pfaffian cubics in $\mathfrak{P}$. For $X$ general the restriction $F\mapsto F_{|X}$ defines a rational map $\rho :\mathpzc{Pf}\dasharrow \mathcal{M}_X$.

\begin{prop}
$\rho $ is birational.
\end{prop}
\pr We have already seen that $\rho $ is dominant. Suppose that two pfaffian cubics $V_1$ and $V_2$, defined by two skew-symmetric  matrices $M_1$ and $M_2$, give the same vector bundle $E$. Restricting the exact sequence (\ref{M}) to $Q$ gives an exact sequence $0\rightarrow \O_{Q}(-1)^6\xrightarrow{\ M_i\ }\O_{Q}^6\xrightarrow{\ p_i\ } E\rightarrow 0$ for $i=1,2$. Since $h^0(E)=h^0(F)=6$, the maps $p_i$ induce an isomorphism on global sections; thus we have a commutative diagram
\[\xymatrix@M=5pt{ 0\ar[r] & \O_{Q}(-1)^6\ar[r]^>>>>{M_1}\ar[d]^A & \O_{Q}^6 \ar[r]^{p_1}\ar[d]^B & E \ar[r]\ar@{=}[d] &  0\\
0 \ar[r]& \O_{Q}(-1)^6\ar[r]^>>>>{M_2} & \O_{Q}^6 \ar[r]^{p_2} & E \ar[r] &  0
}\]where $A$ and $B$ are scalar $6\times 6$ matrices. This implies $V_1=V_2$, hence our assertion.\qed

\medskip	
We now go to the surface $S=X\cap H$ and the moduli space $\mathcal{M}_S$. We have already encountered $\mathcal{M}_S$ in \S 5; it carries a rational Lagrangian fibration $h:\mathcal{M}_S\dasharrow \Pi$, where $\Pi$ is the space $(\cong\P^5)$ of cubics in $H$ containing $S$. The following Proposition tells us that the image of $\mathcal{M}_X$ in $\mathcal{M}_S$ is transverse, in a weak sense, to this fibration:
\begin{prop}
The  rational map $h\rond\operatorname{res}: \mathcal{M}_X\dasharrow \Pi $ is generically finite. 
\end{prop}
 \pr As above, let $\mathfrak{P}$ be the space of cubics in $\P^5$ containing $X$.
 The composite map $h\rond \operatorname{res}\rond \rho :  \allowbreak\mathpzc{Pf}\dasharrow \Pi$ is the restriction to $\mathpzc{Pf}$ of the linear map $p: \mathfrak{P}\dasharrow \Pi$  defined by $p(V)=V\cap H$. This map is defined everywhere except at the point $\mathrm{o}$ corresponding to the cubic $Q\cup H$; in other words, it is the linear projection from $o$. Thus the fiber  $p^{-1}(x)$, for $x\in\Pi$, is the line $\lr{\mathrm{o},x}$. For $x$ general this line intersects $\mathpzc{Pf}$ along a finite set, hence the Proposition.\qed 
 
 \medskip	
As in \S 5, we can fix $S$ and vary $X$, and we obtain a large family of Lagrangian subvarieties of $\mathcal{M}_S$. We do not know whether one can find among them a 5-dimensional family which appears as the fibers of a rational Lagrangian fibration $h':\mathcal{M}_S\dasharrow B$.

\end{document}